\documentclass[11pt]{amsart}

\usepackage{geometry}
\usepackage{amsfonts, adjustbox, amssymb, amscd}
\numberwithin{equation}{section}

\usepackage{bm}
\usepackage{verbatim}
\usepackage{mathrsfs}
\usepackage{graphicx}
\usepackage{tikz-cd}
\usepackage{subcaption}
\usepackage{listings}
\usepackage{subfiles}
\usepackage[toc,page]{appendix}
\usepackage{mathtools}
\usepackage{comment}
\usepackage{enumerate}
\usepackage{enumitem}
\usepackage[all]{xy}

\usepackage{graphicx}
\usepackage{appendix}
\usepackage{hyperref}
\hypersetup{
    colorlinks=true,
    citecolor=red,
    linkcolor=blue,
    filecolor=magenta,      
    urlcolor=red,
}
\newcommand{\bb}{\bm{b}}
\newcommand{\Mm}{{\bf{M}}}

\newcommand{\Qq}{\mathbb{Q}}

\newcommand{\Rr}{\mathbb{R}}

\newcommand{\Center}{\operatorname{center}}

\newcommand{\Exc}{\operatorname{Exc}}

\newcommand{\rk}{\operatorname{rank}}

\newcommand{\ninv}{{\operatorname{ninv}}}

\newcommand{\Supp}{\operatorname{Supp}}

\newcommand{\mult}{\operatorname{mult}}

\newcommand{\Ff}{\mathcal{F}}

\newcommand{\Ll}{{\bf{L}}}


\newcounter{parentnumber}

\newtheorem{thm}{Theorem}[section]
\newtheorem{conj}[thm]{Conjecture}

\newtheorem{lem}[thm]{Lemma}
\newtheorem{prop}[thm]{Proposition}

\theoremstyle{definition}
\newtheorem{defn}[thm]{Definition}

\theoremstyle{definition}
\newtheorem{rem}[thm]{Remark}

\newtheorem{nota}[thm]{Notation}

\theoremstyle{definition}

\begin{document}

\title{Flop between algebraically integrable foliations on potentially klt varieties}
\author{Yifei Chen, Jihao Liu, and Yanze Wang}

\date{\today}

\address{Hua Loo-Keng Key Laboratory of Mathematics and Academy of Mathematics and Systems Science, Chinese Academy of Sciences, No. 55 Zhonguancun East Road, Haidian District, Beijing, 100190, China}
\email{yifeichen@amss.ac.cn}

\address{Department of Mathematics, Peking University, No. 5 Yiheyuan Road, Haidian District, Peking 100871, China}
\email{liujihao@math.pku.edu.cn}

\address{Academy of Mathematics and Systems Science, Chinese Academy of Sciences, No. 55 Zhonguancun East Road, Haidian District, Beijing, 100190, China}
\email{wangyanze@amss.ac.cn}

\pagestyle{myheadings}\markboth{\hfill Yifei Chen, Jihao Liu, and Yanze Wang \hfill}{\hfill Flop between algebraically integrable foliations\hfill}

\begin{abstract}
We prove that for any two minimal models of an lc algebraically integrable foliated triple on potentially klt varieties, there exist small birational models that are connected by a sequence of flops. In particular, any two minimal models of lc algebraically integrable foliated triples on $\Qq$-factorial klt varieties are connected by a sequence of flops. We also discuss the connection between minimal models for possibly non-algebraically integrable foliations on threefolds, assuming the minimal model program for generalized foliated quadruples.
\end{abstract}

\maketitle

\tableofcontents

\section{Introduction}\label{sec: Introduction}

We work over the field of complex numbers $\mathbb C$.

\medskip

\noindent\textbf{Flop between minimal models of varieties.} Let $X$ be a smooth projective variety with a pseudo-effective canonical divisor $K_X$. The minimal model program theory predicts that after a $K_X$-MMP, i.e., a sequence of $K_X$-divisorial contractions and flips, one obtains a minimal model $X_{\min}$ of $X$, such that $X_{\min}$ has $\mathbb{Q}$-factorial terminal singularities and $K_{X_{\min}}$ is nef, i.e., $K_{X_{\min}} \cdot C \geq 0$ for any curve $C$ on $X_{\min}$. However, it is well known that for different $K_X$-MMPs, the minimal model $X_{\min}$ may not be unique. Therefore, it is interesting to ask what the connection is between different minimal models of $X$.

Following the proof of the existence of the minimal model program for klt varieties \cite{BCHM10}, Kawamata proved that the induced birational map $X_1 \dashrightarrow X_2$ between any two minimal models $X_1$ and $X_2$ of a $\mathbb{Q}$-factorial terminal projective variety $X$ can be decomposed into a sequence of $K_{X_1}$-flops \cite[Theorem 1]{Kaw08}. The same lines of the proof also hold when $X$ is $\mathbb{Q}$-factorial klt. Later, following the proof of the existence of the minimal model program for lc varieties \cite{Bir12, HX13, Has19}, Hashizume improved this result to the case when $X$ is $\mathbb{Q}$-factorial lc \cite{Has20}, and also considered the connection between minimal models of $X$ when $X$ is not $\mathbb{Q}$-factorial. More precisely, \cite[Theorem 1.1]{Has20} proved that for any lc projective variety $X$ and two $K_X$-MMPs $X\dashrightarrow X_1$, $X \dashrightarrow X_2$, such that $X_1$ and $X_2$ are minimal models of $X$, there exist small birational morphisms $X_1' \rightarrow X_1$, $X_2' \rightarrow X_2$, and the induced birational map $X_1' \dashrightarrow X_2'$ can be decomposed into a sequence of $K_{X_1'}$-flops. Hashizume's result was later further improved to the case of NQC lc generalized pairs \cite{Cha24}, based on the proof of the existence of the minimal model program for such structures \cite{HL23, LX23a, LX23b, Xie22}.

\medskip

\noindent\textbf{Flop between minimal models of foliations.} In recent years, there have been significant developments in the foundation of the minimal model program theory for foliations. For an lc foliation $ \mathcal{F} $ on a klt variety $ X $, the existence of a $ K_{\mathcal{F}} $-minimal model program (cone theorem, contraction theorem, and the existence of flips) has been proven in dimension 2 \cite{McQ08, Bru15}, in dimension 3 when $\rk\Ff=2$ and $ \mathcal{F} $ is $\mathbb{Q}$-factorial F-dlt \cite{Spi20, CS21, SS22}, in dimension 3 when $\rk\Ff=1$ and $ X $ is $\mathbb{Q}$-factorial \cite{CS20}, and when $ \mathcal{F} $ is algebraically integrable \cite{CHLX23, LMX24b, CHLMSSX24}.

With all the aforementioned progress, it is interesting to ask whether two different minimal models of lc foliations on $\mathbb{Q}$-factorial klt varieties are connected by a sequence of flops, similar to the case of minimal models of varieties. \cite[Theorem 1.1]{JV23} provides a positive answer to this question for $\mathbb{Q}$-factorial F-dlt foliations of rank 2 on threefolds.

The goal of this paper is to provide a positive answer to this question for algebraically integrable foliations in any dimension. More precisely, we have the following:

\begin{thm}\label{thm: main}
Let $\Ff$ be an lc algebraically integrable foliation on a $\Qq$-factorial klt projective variety $X$. Let $\phi_1: (X,\Ff)\dashrightarrow (X_1,\Ff_1)$ and $\phi_2: (X,\Ff)\dashrightarrow (X_2,\Ff_2)$ be two $K_{\Ff}$-MMPs, and let $\alpha: X_1\dashrightarrow X_2$ be the induced birational map. 

Assume that $K_{\Ff}$ is pseudo-effective. Then $\alpha$ can be decomposed into a sequence of $K_{\Ff_1}$-flops.
\end{thm}

In practice, instead of considering $\Ff$ on a projective variety, it is more natural to consider foliated triples $(X,\Ff,B)$ associated with a projective morphism  $X\rightarrow U$ between normal quasi-projective varieties. We have the following result:

\begin{thm}\label{thm: main triple}
Let $(X,\Ff,B)/U$ be a $\Qq$-factorial lc algebraically integrable foliated triple such that $X$ is klt. Let $\phi_1: (X,\Ff,B)\dashrightarrow (X_1,\Ff_1,B_1)$ and $\phi_2: (X,\Ff,B)\dashrightarrow (X_2,\Ff_2,B_2)$ be two $(K_{\Ff}+B)$-MMPs$/U$, and let $\alpha: X_1\dashrightarrow X_2$ be the induced birational map$/U$.

Assume that $K_{\Ff}+B$ is pseudo-effective$/U$. Then $\alpha$ can be decomposed into a sequence of $(K_{\Ff_1}+B_1)$-flops$/U$.
\end{thm}

Inspired by \cite[Theorem 1.2]{Has20}, we also would like to establish the connection between minimal models of lc algebraically integrable foliations on non-$\Qq$-factorial varieties. We have the following result:

\begin{thm}\label{thm: main non q factorial}
Let $(X,\Ff,B)/U$ be an lc algebraically integrable foliated triple such that $X$ is potentially klt. Let $\phi_1: (X,\Ff,B)\dashrightarrow (X_1,\Ff_1,B_1)$ and $\phi_2: (X,\Ff,B)\dashrightarrow (X_2,\Ff_2,B_2)$ be two $(K_{\Ff}+B)$-MMPs$/U$. 

Assume that $K_{\Ff}+B$ is pseudo-effective$/U$. Then there exist small $\Qq$-factorial modifications $(X_1',\Ff_1',B_1')\rightarrow (X_1,\Ff_1,B_1)$ and $(X_2',\Ff_2',B_2')\rightarrow (X_2,\Ff_2,B_2)$, such that the induced birational map $\alpha': X_1'\dashrightarrow X_2'$ can be decomposed into a sequence of $(K_{\Ff_1'}+B_1')$-flops.
\end{thm}

\begin{rem}
The same lines of the proof of our main theorems also apply to algebraically integrable generalized foliated quadruples \cite[Definition 1.2]{LLM23}, possibly with the additional ``NQC" condition. Furthermore, assuming the minimal model program for generalized foliated quadruples on $\Qq$-factorial klt varieties in dimension 3, the same lines of the proof of our main theorems will also extend to (possibly non-algebraically integrable) foliations. In particular, this will provide an alternative proof to \cite[Theorem 1.1]{JV23}. See Section \ref{sec: further discussion} for details. 
\end{rem}

\medskip

\noindent\textit{Sketch of the proofs.} For simplicity, let us assume that $U=\{pt\}$. First, we briefly recall Kawamata's original proof \cite[Theorem 1]{Kaw08} on the flop connection between $\Qq$-factorial terminal pairs $(X_1,B_1)$ and $(X_2,B_2)$ associated with the birational map $\alpha: X_1\dashrightarrow X_2$: we take a general ample divisor $L_2\geq 0$ on $X_2$, consider its strict transform $L_1$ on $X_1$, and run a $(K_{X_1}+B_1+tL_1)$-MMP with scaling of an ample divisor for some $0<t\ll 1$. This MMP terminates with $X_2$ and is a sequence of $(K_{X_1}+B_1)$-flops due to three reasons: 
\begin{enumerate} 
\item Since $0<t\ll 1$ and $(X_1,B_1)$ is klt, $(X_1,B_1+tL_1)$ is still klt. 
\item The bounded length of extremal rays ensures that this MMP is $(K_{X_1}+B_1)$-trivial. 
\item Since $L_1$ is big, this MMP terminates with a good minimal model by \cite{BCHM10}. \end{enumerate} 
We want to follow the idea of \cite[Theorem 1]{Kaw08} to prove Theorem \ref{thm: main triple}. With notations as in Theorem \ref{thm: main triple} and considering (1-3) above, we encounter the following three difficulties: 
\begin{itemize} 
\item[(i)] Due to the failure of Bertini-type theorems, it is possible that $(X_2,\Ff_2,B_2+tL_2)$ is not lc for any ample $\Rr$-divisor $L_2\geq 0$ and any $t>0$. \item[(ii)] Even if $(X_2,\Ff_2,B_2+tL_2)$ is lc, we also do not know how to prove that $(X_1,\Ff_1,B_1+tL_1)$ is lc for any $0<t\ll 1$. \item[(iii)] Even if $(X_1,\Ff_1,B_1+tL_1)$ is lc, $L_1$ is only big. Unlike the usual klt pair case, we still do not know whether a $(K_{\Ff_1}+B_1+tL_1)$-MMP with scaling of an ample divisor will terminate with a good minimal model in general. 
\end{itemize}
Of course, the lengths of extremal rays are still bounded, at least for algebraically integrable foliations \cite[Theorem 3.9]{ACSS21}, \cite[Theorem 2.2.1]{CHLX23}.

For issue (i), the natural idea is to consider the generalized foliated quadruple $(X_2,\Ff_2,B_2;t\Ll:=t\overline{L_2})$ instead of $(X_2,\Ff_2,B_2+tL_2)$, which is a standard way to bypass Bertini-type theorems since \cite{LLM23}. For issue (iii), a key observation is that $K_{\Ff_2}+B_2+tL_2$ is ample. Since $\alpha: X_1\dashrightarrow X_2$ neither extracts nor contracts any divisors, $(X_2,\Ff_2,B_2,t\Ll)$ is a good minimal model of $(X_1,\Ff_1,B_1,t\Ll)$. Now, \cite[Theorem 1.11]{LMX24b} ensures the termination of the $(K_{\Ff_1}+B_1+t\Ll_{X_1})$-MMP with scaling of ample divisors, provided that $(X_1,\Ff_1,B_1,t\Ll)$ is lc. Therefore, we are only left to resolve issue (ii); more precisely, we only need to show that $(X_1,\Ff_1,B_1,t\Ll)$ is lc for some $t>0$.

The difficulty in issue (ii) arises from the fact that $(X_2,\Ff_2,B_2)$ has lc centers, and any small perturbation along the generic point of any lc center may destroy the lc structure. A similar difficulty arose in \cite{Has20} when considering flops connecting lc but not klt pairs $(X_1,B_1)$ and $(X_2,B_2)$. One key idea in \cite[Proof of Theorem 1.1]{Has20} is the observation that, over the generic point of any lc center of $(X_1,B_1)$ or $(X_2,B_2)$, the induced birational maps $X_1\dashrightarrow X_2$ and $X_2\dashrightarrow X_1$ are isomorphisms. In other words, there exist large open subsets $X^0_1\subset X_1$ and $X^0_2\subset X_2$ which contain all lc centers of $(X_1,B_1)$ and $(X_2,B_2)$, respectively, and the induced birational map $X^0_1\dashrightarrow X^0_2$ is an isomorphism. In fact, \cite[Section 4]{Has20} provided counterexamples when no such large open subsets $X^0_1$ and $X^0_2$ exist, i.e., when the minimal models are not obtained as the outputs of the MMPs.

Unfortunately, the arguments in \cite[Proof of Theorem 1.1]{Has20} also do not work in our setting. This is because the lc centers of usual pairs form a proper closed subset of the ambient variety, but for any algebraically integrable foliation $\Ff\neq T_X$, the lc centers of $\Ff$ cover the entire variety $X$. Nevertheless, we have a way to resolve this issue: instead of only considering $(X_1,\Ff_1,B_1)$ and $(X_2,\Ff_2,B_2)$ and finding the corresponding open subsets $X^0_1$ and $X^0_2$, we aim to more concretely use the fact that $(X_1,\Ff_1,B_1)$ and $(X_2,\Ff_2,B_2)$ are outputs of MMPs.

More precisely, let $\phi_i: (X,\Ff,B)\dashrightarrow (X_i,\Ff_i,B_i)$ be $(K_{\Ff}+B)$-MMPs. Then we may take a high foliated log resolution $X'\rightarrow X$ of $(X,\Ff,B)$ so that the induced birational morphisms $X'\rightarrow X_i$ are also foliated log resolutions of $(X_i,\Ff_i,B_i)$. The MMPs $X\dashrightarrow X_i$ lift to MMPs $X'\dashrightarrow X_i'$ with induced $\Qq$-factorial ACSS modifications $(X_i',\Ff_i,B_i')\rightarrow (X_i,\Ff_i,B_i)$. Here are the key observations: \begin{itemize} \item $X_2'\rightarrow X_2$ extracts only finitely many divisors $E_j$ that are lc places of $(X_2,\Ff_2,B_2)$. We may choose $L_2$ so that it does not contain the generic point of the centers of any $E_j$.
\item $E_j$ are also the divisors extracted by $X_1'\rightarrow X_1$, and by a simple discrepancy comparison, $\alpha$ is an isomorphism near the generic point of any $\Center_{X_1}E_j$. Thus, $L_1$ does not contain the generic point of any $\Center_{X_1}E_j$. This implies that $\Ll_{X_1'}$ is the pullback of $L_1$.
\end{itemize}
Since $(X',\Ff',B',t\Ll)$ is lc for any $t>0$, $(X_1',\Ff_1',B_1',t\Ll)$ is lc for any $0<t\ll 1$. Since $\Ll_{X_1'}$ is the pullback of $L_1$, $(X_1,\Ff_1,B_1,t\Ll)$ is lc. This resolves issue (ii), and the proof of Theorem \ref{thm: main triple} follows. Theorems \ref{thm: main} is an easy consequences of Theorem \ref{thm: main triple}, while Theorem \ref{thm: main non q factorial} is a consequence of  Theorem \ref{thm: main triple} together with MMPs lifting to small $\Qq$-factorial modifications (Proposition \ref{prop: small lift mmp}).

\medskip

\noindent\textbf{Acknowledgement}. The authors would like to thank Lingyao Xie for useful discussions. The second author is supported by a start-up funding of Peking University. The first and third authors are partially supported by the NSFC grant No. 12271384.

\section{Preliminaries}\label{sec: preliminaries}

We will adopt the standard notation and definitions on MMP in \cite{KM98,BCHM10} and use them freely. For foliations, foliated triples, generalized foliated quadruples, we adopt the notation and definitions in \cite{LLM23,CHLX23} which generally align with \cite{CS20, ACSS21, CS21} (for foliations and foliated triples) and \cite{BZ16,HL23} (for generalized pairs and $\bb$-divisors), possibly with minor differences. 

\subsection{Foliations}

\begin{defn}[Foliations, {cf. \cite{ACSS21,CS21}}]\label{defn: foliation}
Let $X$ be a normal variety. A \emph{foliation} on $X$ is a coherent sheaf $\Ff\subset T_X$ such that
\begin{enumerate}
    \item $\Ff$ is saturated in $T_X$, i.e. $T_X/\Ff$ is torsion free, and
    \item $\Ff$ is closed under the Lie bracket.
\end{enumerate}

The \emph{rank} of a foliation $\Ff$ on a variety $X$ is the rank of $\Ff$ as a sheaf and is denoted by $\rk\Ff$. 
The \emph{co-rank} of $\Ff$ is $\dim X-\rk\Ff$. The \emph{canonical divisor} of $\Ff$ is a divisor $K_\Ff$ such that $\mathcal{O}_X(-K_{\mathcal{F}})\cong\mathrm{det}(\Ff)$. If $\Ff=0$, then we say that $\Ff$ is a \emph{foliation by points}.

Given any dominant map 
$h: Y\dashrightarrow X$ and a foliation $\mathcal F$ on $X$, we denote by $h^{-1}\Ff$ the \emph{pullback} of $\Ff$ on $Y$ as constructed in \cite[3.2]{Dru21} and say that $h^{-1}\Ff$ is \emph{induced by} $\Ff$. Given any birational map $g: X\dashrightarrow X'$, we denote by $g_\ast \Ff:=(g^{-1})^{-1}\Ff$ the \emph{pushforward} of $\Ff$ on $X'$ and also say that $g_\ast \Ff$ is \emph{induced by} $\Ff$. We say that $\Ff$ is an \emph{algebraically integrable foliation} if there exists a dominant map $f: X\dashrightarrow Z$ such that $\Ff=f^{-1}\Ff_Z$, where $\Ff_Z$ is the foliation by points on $Z$, and we say that $\Ff$ is \emph{induced by} $f$.

A subvariety $S\subset X$ is called \emph{$\Ff$-invariant} if for any open subset $U\subset X$ and any section $\partial\in H^0(U,\Ff)$, we have $\partial(\mathcal{I}_{S\cap U})\subset \mathcal{I}_{S\cap U}$,  where $\mathcal{I}_{S\cap U}$ denotes the ideal sheaf of $S\cap U$ in $U$.  
For any prime divisor $P$ on $X$, we define $\epsilon_{\Ff}(P):=1$ if $P$ is not $\Ff$-invariant and $\epsilon_{\Ff}(P):=0$ if $P$ is $\Ff$-invariant. For any prime divisor $E$ over $X$, we define $\epsilon_{\Ff}(E):=\epsilon_{\Ff_Y}(E)$ where $h: Y\dashrightarrow X$ is a birational map such that $E$ is on $Y$ and $\Ff_Y:=h^{-1}\Ff$.

Suppose that the foliation structure $\Ff$ on $X$ is clear in the context. Then, given an $\mathbb R$-divisor $D = \sum_{i = 1}^k a_iD_i$ where each $D_i$ is a prime Weil divisor,
we denote by $D^{{\rm ninv}} \coloneqq \sum \epsilon_{\mathcal F}(D_i)a_iD_i$ the \emph{non-$\Ff$-invariant part} of $D$ and $D^{{\rm inv}} \coloneqq D-D^{{\rm ninv}}$ the \emph{$\Ff$-invariant part} of $D$.
\end{defn}

\begin{defn}[Singularities]\label{defn: foliation singularity}
Let $(X,\Ff,B)$ be a foliated triple. For any prime divisor $E$ over $X$, let $f: Y\rightarrow X$ be a birational morphism such that $E$ is on $Y$, and suppose that
$$K_{\Ff_Y}+B_Y:=f^*(K_\Ff+B)$$
where $\Ff_Y:=f^{-1}\Ff$. We define $a(E,\Ff,B):=-\mult_EB_Y$ to be the \emph{discrepancy} of $E$ with respect to $(X,\Ff,B)$. We say that $(X,\Ff,B)$ is \emph{lc} if $a(E,\Ff,B)\geq -\epsilon_{\Ff}(E)$ for any prime divisor $E$ over $X$. An \emph{lc place} of $(X,\Ff,B)$ is a prime divisor $E$ over $X$ such that $a(E,\Ff,B)=-\epsilon_{\Ff}(E)$. An \emph{lc center} of $(X,\Ff,B)$ is the center of an lc place of $(X,\Ff,B)$ on $X$.
\end{defn}

\begin{defn}[Potentially klt]\label{defn: potentially klt}
Let $X$ be a normal quasi-projective variety. We say that $X$ is \emph{potentially klt} if $(X,\Delta)$ is klt for some $\Rr$-divisor $\Delta\geq 0$. 
\end{defn}

\subsection{Foliated log smooth}

\begin{defn}[{cf. \cite[3.2 Log canonical foliated pairs]{ACSS21}, \cite[Definition 6.2.1]{CHLX23}}]\label{defn: foliated log smooth}
Let $(X,\Ff,B)/U$ be an algebraically integrable foliated triple. We say that $(X,\Ff,B)$ is \emph{foliated log smooth} if there exists a contraction $f: X\rightarrow Z$ satisfying the following.
\begin{enumerate}
  \item $X$ has at most quotient toric singularities.
  \item $\Ff$ is induced by $f$.
  \item $(X,\Sigma_X)$ is toroidal for some reduced divisor $\Sigma_X$ such that $\Supp B\subset\Sigma_X$.  In particular, $(X,\Supp B)$ is toroidal, and $X$ is $\Qq$-factorial klt.
  \item There exists a log smooth pair $(Z,\Sigma_Z)$ such that $$f: (X,\Sigma_X)\rightarrow (Z,\Sigma_Z)$$ is an equidimensional toroidal contraction.
\end{enumerate}

For any algebraically integrable foliated triple $(X,\Ff,B)$, a \emph{foliated log resolution} of $(X,\Ff,B)$ is a birational morphism $h: X'\rightarrow X$ such that 
$$(X',\Ff':=h^{-1}\Ff,B':=h^{-1}_\ast B+\Exc(h))$$ 
is foliated log smooth, whose existence is guaranteed by \cite{AK00,ACSS21,CHLX23}.
\end{defn}

\begin{defn}[Foliated log smooth model,{ \cite[Definition 4.10]{LMX24b}}]\label{defn: log smooth models}
Let $(X,\Ff,B)$ be an lc algebraically integrable foliated triple and $h: X'\rightarrow X$ a foliated log resolution of $(X,\Ff,B)$. Let $\Ff':=h^{-1}\Ff$, and let $B'\geq 0$ and $E\geq 0$ be two $\Rr$-divisors on $W$ satisfying the following.
\begin{enumerate}
    \item $K_{\Ff'}+B'=h^*(K_\Ff+B)+E$.
    \item $(X',\Ff',B')$ is foliated log smooth and lc.
    \item $E$ is $h$-exceptional.
    \item For any $h$-exceptional prime divisor $D$ such that $$a(D,X,B)>-\epsilon_{\Ff}(E),$$ $D$ is a component of $E$.
\end{enumerate}
We say that $(X',\Ff',B')$ is a \emph{foliated log smooth model} of $(X,\Ff,B)$. 
\end{defn}

\subsection{Models}

\begin{defn}[Minimal models]\label{defn: minimal model}
Let $(X,\Ff,B)/U$ be a foliated triple and $\phi: X\dashrightarrow X'$ a birational map$/U$ which does not extract any divisor. Let $\Ff':=\phi_*\Ff$ and $B':=\phi_*B$. 

We say that $(X',\Ff',B')/U$ is a \emph{weak lc model} (resp. \emph{minimal model}) of $(X,\Ff,B)/U$ if $K_{\Ff'}+B'$ is nef$/U$, and for any prime divisor $D$ on $X$ which is exceptional over $X'$, $a(D,\Ff,B)\leq a(D,\Ff',B')$ (resp. $a(D,\Ff,B)<a(D,\Ff',B')$). We say that $(X',\Ff',B')/U$ is a \emph{semi-ample model} (resp. \emph{good minimal model}) of $(X,\Ff,B)/U$ if $(X',\Ff',B')/U$ is a weak lc model (resp. minimal model) of $(X,\Ff,B)/U$ and $K_{\Ff'}+B'$ is semi-ample$/U$.
\end{defn}

\subsection{Relative Nakayama-Zariski decomposition}

\begin{defn}
    Let $X\rightarrow U$ be a projective morphism between normal quasi-projective varieties and $D$ an $\Rr$-Cartier $\Rr$-divisor on $X$. A birational map$/U$ $\phi: X\dashrightarrow X'$ is called a $D$-MMP$/U$ if $\phi$ is a sequence of steps of a $D$-MMP$/U$, and either there exists a $\phi_*D$-Mori fiber space$/U$ $X'\rightarrow Z$, or $\phi_*D$ is nef$/U$.
\end{defn}

\begin{defn}
    Let $\pi: X\rightarrow U$ be a projective morphism from a normal quasi-projective variety to a quasi-projective variety, $D$ a pseudo-effective$/U$ $\Rr$-Cartier $\Rr$-divisor on $X$, and $P$ a prime divisor on $X$. We define $\sigma_{P}(X/U,D)$ as in \cite[Definition 3.1]{LX23a} by considering $\sigma_{P}(X/U,D)$ as a number in  $[0,+\infty)\cup\{+\infty\}$. We define $N_{\sigma}(X/U,D)=\sum_Q\sigma_Q(X/U,D)Q$
    where the sum runs through all prime divisors on $X$ and consider it as a formal sum of divisors with coefficients in $[0,+\infty)\cup\{+\infty\}$. We say that $D$ is \emph{movable$/U$} if $N_{\sigma}(X/U,D)=0$.
\end{defn}

\begin{lem}\label{lem: output of MMP isomorphic in codimension 1}
  Let $X\rightarrow U$ be a projective morphism between normal quasi-projective varieties and let $D$ be a pseudo-effective$/U$ $\Rr$-Cartier $\Rr$-divisor on $X$. Let $\phi_i: X\dashrightarrow X_i$, $i\in\{1,2\}$ be two $D$-MMP$/U$ and let $D_i:=(\phi_i)_*D$ for $i\in\{1,2\}$. Then:
  \begin{enumerate}
      \item $\phi_2\circ\phi^{-1}_1$ is an isomorphism in codimension $1$.
      \item $D_1$ and $D_2$ are crepant, i.e. for any birational morphisms $p_i: W\rightarrow X_i$, $i\in\{1,2\}$, we have $p_1^*D_1=p_2^*D_2$.
  \end{enumerate}
\end{lem}
\begin{proof}
Since $D$ be a pseudo-effective$/U$, $D_1,D_2$ are nef$/U$, hence $D_1,D_2$ are movable$/U$. By \cite[Lemma 2.25]{LMX24b}, the divisors contracted by $\phi_i$ are exactly $\Supp N_{\sigma}(X/U,D)$. This implies (1). By (1), $p_1^*D_1-p_2^*D_2$ is nef$/X_2$ and exceptional$/X_2$, and $p_2^*D_2-p_1^*D_1$ is nef$/X_1$ and exceptional$/X_1$. By applying the negativity lemma twice, $p_1^*D_1=p_2^*D_2$. This implies (2).
\end{proof}

\subsection{Generalized pairs and generalized foliated quadruples}

\begin{rem}
For $\bb$-divisors and generalized pairs, we will follow the notations and definitions in \cite{BZ16, HL23}. For generalized foliated quadruples, we will follow \cite{CHLX23}.

Generalized pairs and generalized foliated quadruples are very technical concepts. To make the statements in this paper more concise, for most of the results we cite and prove, we refer only to the foliated triple version and briefly mention the detailed versions in Section \ref{sec: further discussion}. We refer the reader to \cite[Appendix A]{LMX24b} and \cite{CHLX23} for results on generalized foliated quadruples.

The concepts defined above (singularities of foliations, foliated log smooth, foliated log smooth model, minimal models, etc.) can be similarly defined for generalized foliated quadruples. For the reader's convenience, we omit those technical definitions.
\end{rem}

\begin{defn}[NQC]
    Let $X\rightarrow U$ be a projective morphism from a normal quasi-projective variety to a variety. Let $D$ be a nef $\Rr$-Cartier $\Rr$-divisor on $X$ and $\Mm$ a nef $\bb$-divisor on $X$.

    We say that $D$ is NQC$/U$ if $D=\sum d_iD_i$, where each $d_i\geq 0$ and each $D_i$ is a nef$/U$ Cartier divisor. We say that $\Mm$ is NQC$/U$ if $\Mm=\sum \mu_i\Mm_i$, where each $\mu_i\geq 0$ and each $\Mm_i$ is a nef$/U$ Cartier $\bb$-divisor.
\end{defn}

\section{Proof of Main Theorem}\label{sec: proof}

\begin{proof}[Proof of Theorem \ref{thm: main triple}]
By Lemma \ref{lem: output of MMP isomorphic in codimension 1}, $\alpha$ is an isomorphism in codimension $1$ and $\alpha_*B_1=B_2$.

Let $h: X'\rightarrow X$ and $h_i: X'\rightarrow X_i$, $i\in\{1,2\}$ be birational morphisms such that $h_i$ is a foliated log resolution of $(X_i,\Ff_i,B_i)$ and $h$ is a foliated log resolution of $(X,\Ff,B)$. Let $B':=h^{-1}_*B+\Exc(h)^{\ninv}$. Since $(X,\Ff,B)$ is lc, $(X',\Ff',B')$ is a foliated log smooth model of $(X,\Ff,B)$ and $B'=(h_i^{-1})_*B_i+F_i$ for some $F_i\geq 0$ that are exceptional$/X_i$. 

By \cite[Lemma 4.13]{LMX24b}, for $i\in\{1,2\}$, we may run a $(K_{\Ff'}+B')$-MMP$/X_i$ with scaling of an ample divisor, which terminates with a good minimal model $(X_i',\Ff_i',B_i')/X_i$ of $(X',\Ff',B')/X_i$ with induced birational morphism $g_i: X_i'\rightarrow X_i$ and birational map $\phi_i': X'\dashrightarrow X_i'$, such that $K_{\Ff_i'}+B_i'=g_i^*(K_{\Ff_i}+B_i)$. Since $K_{\Ff_i}+B_i$ is nef$/U$, $K_{\Ff_i'}+B_i'$ is nef$/U$. Thus $(X_i',\Ff_i',B_i')/U$ is a minimal model of $(X',\Ff',B')/U$. 

Let $\alpha': X_1'\dashrightarrow X_2'$ be the induced birational map$/U$. By Lemma \ref{lem: output of MMP isomorphic in codimension 1}, $\alpha'$ is an isomorphism in codimension $1$. Therefore, the divisors extracted by $g_1$ are exactly the divisors extracted by $g_2$. Let $E_1,\dots,E_m$ be these divisors. Then each $E_j$ is an lc place of $(X_i,\Ff_i,B_i)$ for $i\in\{1,2\}$, so each $E_j$ is also an lc place of $(X,\Ff,B)$ as
$$-\epsilon_{\Ff}(E_j)\leq a(E_j,\Ff,B)\leq a(E_j,\Ff_1,B_1)=-\epsilon_{\Ff}(E_j).$$
Therefore, $\phi_i$ is an isomorphism near $\Center_XE_j$ for $i\in\{1,2\}$ and any $j$, so $\alpha$ (resp. $\alpha^{-1}$) is an isomorphism near $\Center_{X_1}E_j$ (resp. $\Center_{X_2}E_j$) for each $j$.

Let $L_2\geq 0$ be an ample divisor on $X_2$ such that $L_2$ does not contain the generic point of $\Center_{X_2}E_j$ for any $j$, and $\Ll:=\overline{L_2}$. Since $L$ is ample and $K_{\Ff_2}+B_2$ is nef$/U$, for any $s>0$, $X_2$ is the ample model$/U$ of $K_{\Ff_2}+B_2+s\Ll_{X_2}$. 

Since $(X',\Ff',B')/U$ is lc and $\Ll$ descends to $X'$, $(X',\Ff',B',s\Ll)/U$ is lc (as a generalized foliated quadruple) for any $s>0$. For any $0<s\ll 1$, $\phi_1'$ is a sequence of steps of a $(K_{\Ff'}+B'+s\Ll_{X'})$-MMP$/U$, hence $(X_1',\Ff_1',B_1',s\Ll)/U$ is lc. Since $\alpha$ (resp. $\alpha^{-1}$) is an isomorphism near $\Center_{X_1}E_j$ (resp. $\Center_{X_2}E_j$) for each $j$, $\Ll_{X_1}=\alpha^{-1}_{\ast}L_2$ does not contain the generic point of $\Center_{X_1}E_j$ for any $j$. Since $X$ is $\Qq$-factorial, $X_1$ is $\Qq$-factorial. Thus
$$\Ll_{X_1'}=g_1^*\Ll_{X_1}.$$
Since $K_{\Ff_1'}+B_1'=g_1^*(K_{\Ff_1}+B_1)$, we have
$$K_{\Ff_1'}+B_1'+s\Ll_{X_1'}=g_1^*(K_{\Ff_1}+B_1+s\Ll_{X_1}),$$
hence $(X_1,\Ff_1,B_1,s\Ll)/U$ is lc for any $0<s\ll 1$. We let $s_0<1$ be a positive real number such that  $(X_1,\Ff_1,B_1,s\Ll)/U$ is lc for any $0<s\leq s_0$.

Since $h_1$ is a foliated log resolution of $(X_1,\Ff_1,B_1)$ and $\Ll$ descends to $X'$, $h_1$ is a foliated log resolution of $(X_1,\Ff_1,B_1,s\Ll)$, and $(X',\Ff',B',s\Ll)$ is a foliated log smooth model of $(X_1,\Ff_1,B_1,s\Ll)$ for any $0<s\leq s_0$. By our construction, $K_{\Ff_2}+B_2+s\Ll_{X_2}$ is the ample model$/U$ of $K_{\Ff'}+B'+s\Ll_{X'}$ for any $0<s\leq s_0$. Therefore, $(X_2,\Ff_2,B_2,s\Ll)/U$ is a semi-ample model of $(X',\Ff',B',s\Ll)/U$ for any $0<s\leq s_0$. By \cite[Lemma A.28]{LMX24b} and since $\alpha$ does not extract any divisor, $(X_2,\Ff_2,B_2,s\Ll)/U$ is a semi-ample model of $(X_1,\Ff_1,B_1,s\Ll)/U$ for any $0<s\leq s_0$. Since $\alpha$ does not contract any divisor, $(X_2,\Ff_2,B_2,s\Ll)/U$ is a good minimal model of $(X_1,\Ff_1,B_1,s\Ll)/U$ for any $0<s\leq s_0$. In particular,  $(X_1,\Ff_1,B_1,s\Ll)/U$ has a good minimal model for any $0<s\leq s_0$. 

By \cite[Theorem 7.2]{LMX24b}, $X_1$ is $\Qq$-factorial klt. By \cite[Theorem A.13]{LMX24b}, for any $0<s\leq s_0$, we may run a $(K_{\Ff_1}+B_1+s\Ll_{X_1})$-MMP$/U$ which terminates with a good minimal model $(X_s,\Ff_s,B_s,s\Ll)/U$ of $(X_1,\Ff_1,B_1,s\Ll)/U$, such that $X_s$ is $\Qq$-factorial klt. Since $K_{\Ff_2}+B_2+s\Ll_{X_2}$ is ample$/U$, by \cite[Lemma A.25]{LMX24b}, there exists an induced morphism $\psi_s: X_s\rightarrow X_2$. Since $\alpha$ is small and the induced birational map $\alpha_s: X_1\dashrightarrow X_s$ does not extract any divisor, $\psi_s$ and $\alpha_s$ are small. However, since $X_s$ is also $\Qq$-factorial klt, $\psi_s$ is the identity morphism, $\alpha_s=\alpha$, and $X_s=X_2$.

By \cite[Theorem 1.12]{LMX24b}, $K_{\Ff_1}+B_1$ is NQC$/U$. Since $\alpha_s=\alpha$ can be decomposed into a sequence of steps of a
$$\left((K_{\Ff_1}+B_1+s_0\Ll_{X_1})+\left(\frac{s_0}{s}-1\right)(K_{\Ff_1}+B_1)\right)\text{-MMP}/U$$
for any $0<s\leq s_0$. By \cite[Lemma B.6]{LMX24b}, for $0 <s \gg s_0$, $\alpha=\alpha_s$ can be decomposed into a sequence of steps of a $(K_{\Ff_1}+B_1+s\Ll_{X_1})$-MMP$/U$, and each step is $(K_{\Ff_1}+B_1)$-trivial. In particular, $\alpha=\alpha_s$ can be decomposed into a sequence of $(K_{\Ff_1}+B_1)$-flops$/U$.
\end{proof}

\begin{proof}[Proof of Theorem \ref{thm: main}]
    It is a special case of Theorem \ref{thm: main triple}.
\end{proof}

To prove Theorem \ref{thm: main non q factorial}, we need the following proposition, Proposition \ref{prop: small lift mmp}, on MMP lifting. The difference between  Proposition \ref{prop: small lift mmp} and \cite[Proposition 8.2]{LMX24b} is that \cite[Proposition 8.2]{LMX24b} lifts the MMP to $\Qq$-factorial ACSS models, while  Proposition \ref{prop: small lift mmp} lifts the MMP to small $\Qq$-factorial modifications. 

\begin{prop}\label{prop: small lift mmp}
    Let $(X,\Ff,B)/U$ be an lc algebraically integrable foliated triple. Let $\mathcal{P}:$
$$(X,\Ff,B):=(X_0,\Ff_0,B_0)\dashrightarrow (X_1,\Ff_1,B_1)\dashrightarrow\dots\dashrightarrow (X_n,\Ff_n,B_n)\dashrightarrow\dots$$
be a (possibly infinite) sequence of $(K_{\Ff}+B)$-MMP$/U$. For each $i\geq 0$, we let $\psi_i: X_i\rightarrow T_i$ and $\psi_i^+:X_{i+1}\rightarrow T_{i}$ be the $(i+1)$-th step of this MMP and let $\phi_i:=(\psi_{i}^+)^{-1}\circ\psi_i: X_i\dashrightarrow X_{i+1}$ be the induced birational map. Let $A$ be an ample$/U$ $\Rr$-divisor on $X$ and let $A_i$ be the image of $A$ on $X_i$ for each $i$.

Assume that $X$ is potentially klt. Let $h: Y\rightarrow X$ be a small $\Qq$-factorial modifications of $X$, $\Ff_Y:=h^{-1}\Ff$, and $B_Y:=h^{-1}_*B$. Then there exist a (possibly infinite) sequence $\mathcal{P}_Y$ of birational maps 
$$(Y,\Ff_Y,B_Y):=(Y_0,\Ff_{Y_0},B_{Y_0})\dashrightarrow (Y_1,\Ff_{Y_1},B_{Y_1})\dashrightarrow\dots\dashrightarrow (Y_n,\Ff_{Y_n},B_{Y_n})\dashrightarrow\dots$$
satisfying the following. Let $\phi_{i,Y}: Y_i\dashrightarrow Y_{i+1}$ be the induced birational map. Then:
\begin{enumerate}
\item For any $i\geq 0$, there exist a small $\Qq$-factorial modifications $h_i: Y_i\rightarrow X_i$ such that $\Ff_{Y_i}=h_i^{-1}\Ff_i$, $B_{Y_i}=(h_i^{-1})_*B_i$, and $h_0=h$.
\item For any $i\geq 0$, $h_{i+1}\circ\phi_{i,Y}=\phi_i\circ h_i$.
\item For any $i\geq 0$, $\phi_{i,Y}$ is a $(K_{\Ff_i}+B_{Y_i})$-MMP$/T_i$ and $(Y_{i+1},\Ff_{Y_{i+1}},B_{Y_{i+1}})/T_i$ is the output of this MMP, such that $\phi_{i,Y}$ is not the identity map.
\item $\mathcal{P}_Y$ is a sequence of steps of a $(K_{\Ff_Y}+B_Y)$-MMP$/U$.
\item Suppose that $\mathcal{P}$ is an MMP$/U$ with scaling of $A$. Let $A_{Y}:=h^*A$ and let $A_{Y_i}$ the image of $A_Y$ on $Y_i$ for each $i$. Let
$$\lambda_i:=\inf\{t\geq 0\mid K_{\Ff_i}+B_i+tA_i\text{ is nef}/U\}$$
be the $(i+1)$-th scaling number. Then:
\begin{enumerate}
    \item $\phi_{i,Y}$ is a sequence of steps of a $(K_{\Ff_{Y_i}}+B_{Y_i})$-MMP$/U$ with scaling of $A_{Y_i}$, and the scaling number of each step of $\phi_{i,Y}$ is $\lambda_i$.
    \item $\mathcal{P}_Y$ is sequence of steps of a $(K_{\Ff_Y}+B_Y)$-MMP$/U$ with scaling of $A_Y$.
\end{enumerate}
\end{enumerate}
\end{prop}
\begin{proof}
Since (4) follows from (3) and (5.b) follows from (5.a), we only need to prove (1)(2)(3) and (5.a).

Let $n$ be a non-negative integer. We prove the proposition by induction on $n$ and under the induction hypothesis that we have already constructed $(Y_i,\Ff_i,B_i)/U$ and $h_i$ for any $i\leq n$ and $\phi_{i,Y}$ for any $i\leq n-1$ which satisfy (1)(2)(3)(5). When $n=0$, this follows from our assumption, so we may assume that $n>0$. We need to construct $\phi_{n,Y},h_{n+1}$, and $(Y_{n+1},\Ff_{n+1},B_{n+1})/U$.

We let $H_n$ be a supporting function of the extremal ray$/U$ contracted by $\psi_n$ and let  
$$L_n:=H_n-(K_{\Ff_n}+B_n),$$
such that $L_n=\lambda_nA_n$ if $\mathcal{P}$ is an MMP$/U$ with scaling of $A$. Then $L_n$ is ample$/T_n$. Since $K_{\Ff_{Y_n}}+B_{Y_n}+h_n^*L_n\equiv_{T_n}0$ and $K_{\Ff_{Y_n}}+B_{Y_n}+h_n^*L_n$ is nef$/U$, we may run a $(K_{\Ff_{Y_n}}+B_{Y_n})$-MMP$/T_n$ with scaling of an ample divisor, which is also a $(K_{\Ff_{Y_n}}+B_{Y_n})$-MMP$/T_n$ with scaling of $h_n^*L_n$, and is also a sequence of steps of a $(K_{\Ff_{Y_n}}+B_{Y_n})$-MMP$/U$ with scaling of $h_n^*L_n$. Since $(X_{n+1},\Ff_{n+1},B_{n+1})/T_n$ is a semi-ample model of $(Y_{n},\Ff_{Y_n},B_{Y_n})/T_n$, by \cite[Theorem 1.11]{LMX24b}, the  $(K_{\Ff_{Y_n}}+B_{Y_n})$-MMP$/T_n$ with scaling of an ample divisor terminates with a good minimal model $(Y_{n+1},\Ff_{Y_{n+1}},B_{Y_{n+1}})/T_n$ of $(Y_{n},\Ff_{Y_{n}},B_{Y_{n}})/T_n$. Since $X_{n+1}$ is the ample model$/T_n$ of $K_{\Ff_n}+B_n$, $X_{n+1}$ is also the ample model$/T_n$ of $K_{\Ff_{Y_{n+1}}}+B_{Y_{n+1}}$, so there exists an induced birational morphism $h_{n+1}: Y_{n+1}\rightarrow X_{n+1}$. Since $K_{Y_n}+B_{Y_n}$ is not nef$/T_n$ and $K_{Y_{n+1}}+B_{Y_{n+1}}$ is nef$/T_n$, $\phi_{i,Y}$ is not the identity map.

Let $\phi_{n,Y}: Y_n\dashrightarrow Y_{n+1}$ be the induced birational map. By \cite[Lemma 2.25]{LMX24b}, the divisors contracted by $\phi_n$ are exactly $\Supp N_{\sigma}(X_n/T_n,K_{\Ff_n}+B_n)$ and the divisors contracted by $\phi_{n,Y}$ are exactly $\Supp N_{\sigma}(Y_n/T_n,K_{\Ff_{Y_n}}+B_{Y_n})$. By \cite[Lemma 3.4(2)(3)]{LX23a}, 
$$\Supp N_{\sigma}(Y_n/T_n,K_{\Ff_{Y_n}}+B_{Y_n})=\Supp N_{\sigma}(X_n/T_n,K_{\Ff_n}+B_n).$$
Thus $h_{n+1}$ is small, which implies (1) for $n+1$. (2)(3) for $n+1$ follow immediately from our construction. Since $H_n\sim_{\mathbb R,T_n}0$, $\phi_{n,Y}$ is $(h_n^*H_n)$-trivial, so (5.a) for $n+1$ immediately follows. Thus (1)(2)(3) and (5.a) follow from induction on $n$ and the proposition follows.
\end{proof}

\begin{proof}[Proof of Theorem \ref{thm: main non q factorial}]
Let $h: X'\rightarrow X$ be a small $\Qq$-factorial modification, $\Ff':=h^{-1}\Ff$, and $B':=h^{-1}_*B$. By Proposition \ref{prop: small lift mmp}, for $i\in\{1,2\}$, there exists a $(K_{\Ff'}+B')$-MMP$/U$ $\phi_i': (X',\Ff',B')\dashrightarrow (X_i',\Ff_i',B_i')$ with induced small $\Qq$-factorialization $h_i: X_i'\rightarrow X_i$, such that $\Ff_i'=h_i^{-1}\Ff_i$ and $B_i'=(h_i^{-1})_*B_i$. By Theorem \ref{thm: main triple}, the induced birational map $\alpha': X_1'\dashrightarrow X_2'$ can be decomposed into a sequence of $(K_{\Ff_1'}+B_1')$-flops$/U$. The theorem follows.
\end{proof}

\section{Further discussions}\label{sec: further discussion}

\subsection{Generalized foliated quadruple versions}

The main theorems of our paper also hold for NQC lc algebraically integrable generalized foliated quadruples on $\Qq$-factorial klt varieties. We cannot remove the NQC condition because the proof of Theorem \ref{thm: main triple} relies on Shokurov polytopes \cite[Theorems 1.12]{LMX24b}, which also apply to NQC lc algebraically integrable generalized foliated quadruples \cite[Theorem A.14]{LMX24b}, but no longer hold in the non-NQC case.

\begin{defn}
A \emph{generalized foliated quadruple} $(X,\Ff,B,\Mm)/U$ consists of a normal quasi-projective variety $X$, a foliation $\Ff$ on $X$, an $\Rr$-divisor $B\geq 0$ on $X$, a projective morphism $X\rightarrow U$, and a nef$/U$ $\bb$-divisor $\Mm$, such that $K_{\Ff}+B+\Mm_X$ is $\mathbb R$-Cartier. We say that $(X,\Ff,B,\Mm)/U$ is \emph{NQC} if $\Mm$ is NQC$/U$. 
\end{defn}

\begin{thm}[Theorem \ref{thm: main triple}]\label{thm: main triple gfq}
Let $(X,\Ff,B,\Mm)/U$ be a $\Qq$-factorial NQC lc algebraically integrable generalized foliated quadruple such that $X$ is klt. Let $\phi_1: (X,\Ff,B,\Mm)\dashrightarrow (X_1,\Ff_1,B_1,\Mm)$ and $\phi_2: (X,\Ff,B,\Mm)\dashrightarrow (X_2,\Ff_2,B_2,\Mm)$ be two $(K_{\Ff}+B+\Mm_X)$-MMPs$/U$, and let $\alpha: X_1\dashrightarrow X_2$ be the induced birational map$/U$.

Assume that $K_{\Ff}+B+\Mm_X$ is pseudo-effective$/U$. Then $\alpha$ can be decomposed into a sequence of $(K_{\Ff_1}+B_1+\Mm_{X_1})$-flops$/U$.
\end{thm}
\begin{proof}
    The proof follows from the same lines of the proof of Theorem \ref{thm: main triple} except that we replace \cite[Lemma 4.13]{LMX24b} with \cite[Lemma A.29]{LMX24b} and \cite[Theorem 1.12]{LMX24b} with \cite[Theorem A.14]{LMX24b} respectively.
\end{proof}

\begin{thm}[Theorem \ref{thm: main non q factorial}]\label{thm: main non q factorial gfq}
Let $(X,\Ff,B,\Mm)/U$ be an NQC lc algebraically integrable generalized foliated quadruple such that $X$ is potentially klt. Let $\phi_1: (X,\Ff,B,\Mm)\dashrightarrow (X_1,\Ff_1,B_1,\Mm)$ and $\phi_2: (X,\Ff,B,\Mm)\dashrightarrow (X_2,\Ff_2,B_2,\Mm)$ be two $(K_{\Ff}+B+\Mm_X)$-MMPs$/U$. 

Assume that $K_{\Ff}+B+\Mm_X$ is pseudo-effective$/U$. Then there exist small $\Qq$-factorial modifications $(X_1',\Ff_1',B_1',\Mm)\rightarrow (X_1,\Ff_1,B_1,\Mm)$ and $(X_2',\Ff_2',B_2',\Mm)\rightarrow (X_2,\Ff_2,B_2,\Mm)$, such that the induced birational map $\alpha': X_1'\dashrightarrow X_2'$ can be decomposed into a sequence of $(K_{\Ff_1'}+B_1'+\Mm_{X_1'})$-flops.
\end{thm}
\begin{proof}
    The proof follows from the same lines of the proof of Theorem \ref{thm: main non q factorial} except that we replace Theorem \ref{thm: main triple} with Theorem \ref{thm: main triple gfq} and replace Proposition \ref{prop: small lift mmp} with the corresponding generalized foliated quadruple version. The generalized foliated quadruple version of Proposition \ref{prop: small lift mmp} also from the same lines of the proof, except that we replace \cite[Theorem 1.11]{LMX24b} with \cite[Theorem A.13]{LMX24b}.
\end{proof}

\subsection{Non-algebraically integrable foliations on threefolds}

With the flop connection theorems for algebraically integrable foliations settled, it is natural to ask whether similar theorems hold for non-algebraically integrable foliations. Indeed, we can deduce similar theorems, provided that the minimal model program holds for generalized foliated quadruples. We first consider the following conjecture.

\begin{conj}\label{conj: mmp gfq threefold}
    Let $(X,\Ff,B,\Mm)/U$ be a $\Qq$-factorial NQC lc generalized foliated quadruple such that $\dim X=3$ and $X$ is klt. Let $A$ be an ample$/U$ $\Rr$-divisor. Then:
    \begin{enumerate}
    \item The cone theorem (including the boundedness of length of extremal rays and the finiteness of $(K_{\Ff}+B+\Mm_X+A)$-negative extremal rays$/U$ in $\overline{NE}(X/U)$), contraction theorem, and the existence of flips hold for $(X,\Ff,B,\Mm)/U$.
        \item For any sequence of steps of a $(K_{\Ff}+B+\Mm_X)$-MMP$/U$ $(X,\Ff,B,\Mm)\dashrightarrow (X',\Ff',B',\Mm)$, $X'$ is $\Qq$-factorial klt. 
        \item If $(X,\Ff,B,\Mm)/U$ has a minimal model, then we may run a $(K_{\Ff}+B+\Mm_X)$-MMP$/U$ with scaling of an ample$/U$ $\Rr$-divisor and any such MMP terminates with a minimal model of $(X,\Ff,B,\Mm)/U$.
        \item $(X,\Ff,B+A,\Mm)/U$ has a minimal model.
    \end{enumerate}
\end{conj}

Now we consider the following variations of Conjecture \ref{conj: mmp gfq threefold}:

\begin{nota}\label{nota: special conjecture}
    We use the following subscripts to describe the special cases of Conjecture \ref{conj: mmp gfq threefold}:
    \begin{enumerate}
        \item ``1": Conjecture \ref{conj: mmp gfq threefold} when $\rk\Ff=1$.
        \item ``2": Conjecture \ref{conj: mmp gfq threefold} when $\rk\Ff=2$.
        \item ``Fdlt": Conjecture \ref{conj: mmp gfq threefold} when $\rk\Ff=2$ and $(X,\Ff,B,\Mm)/U$ is F-dlt, with the extra conjecture that any for any sequence of steps of a $(K_{\Ff}+B+\Mm_X)$-MMP$/U$ $(X,\Ff,B,\Mm)\dashrightarrow (X',\Ff',B',\Mm)/U$, $(X',\Ff',B',\Mm)$ is F-dlt.
        \item ``Proj": Conjecture \ref{conj: mmp gfq threefold} when $U=\{pt\}$.
        \item ``$\emptyset$": Conjecture \ref{conj: mmp gfq threefold}.
    \end{enumerate}
    We let $\Lambda$ be the set of the following subscripts: $\emptyset$; $1$; $2$; Fdlt; Proj; 1, Proj; 2, Proj; Fdlt, Proj.
\end{nota}

\begin{thm}\label{thm: assume gfq mmp get flop connection}
Let $\lambda\in\Lambda$ be a subscript. Assume Conjecture \ref{conj: mmp gfq threefold}{$_\lambda$} holds. 

Let $(X,\Ff,B)/U$ be a $\Qq$-factorial lc foliated triple such that $\dim X=3$ and $X$ is klt. Let $\phi_1: (X,\Ff,B)\dashrightarrow (X_1,\Ff_1,B_1)$ and $\phi_2: (X,\Ff,B)\dashrightarrow (X_2,\Ff_2,B_2)$ be two $(K_{\Ff}+B)$-MMPs$/U$, and let $\alpha: X_1\dashrightarrow X_2$ be the induced birational map$/U$. Assume that $K_{\Ff}+B$ is pseudo-effective$/U$. Moreover:
\begin{enumerate}
    \item If $\lambda$ contains $1$ (resp. $2$), then we assume that $\rk\Ff=1$ (resp. $\rk\Ff=2)$.
    \item If $\lambda$ contains Fdlt, then we assume that $\rk\Ff=2$ and $(X,\Ff,B)$ is F-dlt.
    \item If $\lambda$ contains Proj, then we assume that $U=\{pt\}$.
\end{enumerate}
Then $\alpha$ can be decomposed into a sequence of $(K_{\Ff_1}+B_1)$-flops$/U$. 
\end{thm}
\noindent\textit{Sketch of the proof.} The proof follows from the same lines of the proof of Theorem \ref{thm: main triple} with the following major modifications:
\begin{enumerate} 
\item When we take the foliated log resolution $h: X'\rightarrow X$, it shall be replaced with the foliated log resolution in the sense of \cite[Definition 3.1]{CS21} if $\rk\Ff=2$, whose existence is guaranteed by \cite{Can04}, and in the sense of \cite[Definition 4.4]{LLM23} if $\rk\Ff=1$, whose existence is guaranteed by \cite{MP13}. See also \cite[Theorem 4.5]{LLM23} for an explanation. 
\item The results \cite[Lemma 4.13, A.25, A.28]{LMX24b} used in the proof of Theorem \ref{thm: main triple} shall be replaced with the corresponding non-algebraically integrable version in dimension 3. However, all these results essentially only rely on the ``very exceptional minimal model program", which in turn only relies on Birkar's general negativity lemma \cite[Lemma 3.3]{Bir12} (which always holds) and Conjecture \ref{conj: mmp gfq threefold}{$_\lambda$}(3-4). \item The result \cite[Theorem 7.2]{LMX24b} used in the proof of Theorem \ref{thm: main triple} can be replaced with Conjecture \ref{conj: mmp gfq threefold}{$_\lambda$}(2) (or, if $\lambda$ contains Fdlt, with Notation \ref{nota: special conjecture}(3)). 
\item The result \cite[Theorem A.13]{LMX24b} used in the proof of Theorem \ref{thm: main triple} can be replaced with Conjecture \ref{conj: mmp gfq threefold}{$_\lambda$}(3). 
\item The results \cite[Theorem 1.12, Lemma B.6]{LMX24b} used in the proof of Theorem \ref{thm: main triple} essentially only rely on \cite[Lemma 5.3]{HLS19} (which always hold, or one can simply apply \cite[Theorem 1.2]{LMX24a}), the finiteness of $(K_{\Ff}+B+\Mm_X+A)$-negative extremal rays$/U$, and the boundedness of the length of extremal rays. The latter two conditions follow from Conjecture \ref{conj: mmp gfq threefold}{$_\lambda$}(1).\hfill$\Box$ \end{enumerate}

In particular, Theorem \ref{thm: assume gfq mmp get flop connection}(2) and Conjecture \ref{conj: mmp gfq threefold}$_{\text{Fdlt,Proj}}$ together will provide a new proof of \cite[Theorem 1.1]{JV23}. We can also state a generalized foliated quadruple version of Theorem \ref{thm: assume gfq mmp get flop connection}, and the proof will be similar. Again, due to technicalities, we omit the statements.

\noindent\textbf{Postscript Remark.} After we finished the first draft of this paper, we notice that \cite[Theorem 8.1]{CM24} proved the flop connection theorem for lc foliated triples $(X,\Ff,B)/U$ such that $\dim X=3$, $\rk\Ff=2$, $B$ is a $\Qq$-divisor, $(X,B)$ is klt, and $U=\{pt\}$. Their proof relies on the proof of a variation of Conjecture \ref{conj: mmp gfq threefold} when $\rk\Ff=2$, $B$ is a $\Qq$-divisor, $\Mm$ is a $\Qq$-$\bb$-divisor and $\bb$-semi-ample, $(X,B,\Mm)$ is klt, and $U=\{pt\}$. The proof follows from the same logic as the sketch of the proof of Theorem \ref{thm: assume gfq mmp get flop connection}.

\subsection{Final remarks}
\begin{rem} In our main theorems, we require that $X$ is $\Qq$-factorial klt or potentially klt. This is because we only know the existence of the minimal model program for algebraically integrable foliations under these extra conditions \cite[Theorems 1.2, 1.3]{LMX24b}, \cite[Theorem 1.11]{CHLMSSX24}. We expect the ``potentially klt" condition to be removed if the existence of the minimal model program is proven for any lc algebraically integrable foliated triples. \end{rem}

\begin{rem} It seems that we still do not know whether flops will connect minimal models for non-NQC lc generalized pairs, even when the ambient variety is $\Qq$-factorial klt. For non-NQC $\Qq$-factorial klt generalized pairs, although we cannot find a written proof, a positive answer can be deduced following the same lines of the proof in \cite[Proof of Theorem 1]{Kaw08}. \end{rem}

\begin{rem} The dual question to ``flops connecting minimal models" is the Sarkisov program, which provides a way to connect Mori fiber spaces. The Sarkisov program was established in \cite[Theorem 1.3]{HM13} for $\Qq$-factorial klt pairs, in \cite[Theorem 1.5]{Liu21} for $\Qq$-factorial klt generalized pairs, and in \cite[Theorem 1.1]{Mas24} for $\Qq$-factorial F-dlt foliated triples $(X,\Ff,B)/U$ with $\dim X=3$, $\rk\Ff=2$, and $\lfloor B\rfloor=0$.

With the ``flop connection part" settled in this paper, it is interesting to ask whether we can establish the Sarkisov program for lc algebraically integrable foliated triples on $\Qq$-factorial klt varieties. This should follow from the log geography of minimal models, which is yet to be established for algebraically integrable foliations. \end{rem}

\end{document}